\documentclass{amsart}
\usepackage{epsfig,psfrag,amsmath,amssymb,longtable,dcolumn}
\newtheorem{theorem}{Theorem}

\newtheorem{lemma}{Lemma}

\theoremstyle{definition}

\theoremstyle{remark}

\newcommand{\llangle}{\langle\langle}
\newcommand{\rrangle}{\rangle\rangle}

\def\co{\colon\thinspace}

\newcommand{\begriff}[1]{\textbf{#1}}                  % -"-  keinen -"-

\begin{document}

\title[Minimal generating sets of non-modular invariant rings]{Minimal generating sets of non-modular invariant rings of finite groups}

\begin{abstract}
        It is a classical problem to compute a minimal set of invariant polynomials
        generating the invariant ring of a finite group as an algebra. We present here an
        algorithm for the computation of minimal generating sets in the non-modular case. 
        Apart from very few explicit computations of Gr\"obner bases, the algorithm only 
        involves very basic operations.
        
        As a test bed for comparative benchmarks, we use transitive permutation groups on $7$ and 
        $8$ variables. In most examples, our algorithm implemented in \textsc{Singular} works much 
        faster than the one used in \textsc{Magma}, namely by factors between $50$ and $1000$. We 
        also compute some further examples on more than $8$ variables, including a minimal generating
        set for the natural action of the cyclic group of order $11$ in characteristic $0$ and 
        of order up to $15$ in small prime characteristic. 
        
        We also apply our algorithm to the computation of irreducible secondary invariants.

\end{abstract}
\author{Simon A. King}
\address{Simon A. King\\
Mathematisches Forschungsinstitut Oberwolfach\\
Schwarzwaldstr. 9--11\\
D-77709 Oberwolfach\\
Germany}
\email{king@mfo.de}
\maketitle

\section{Introduction}
Let $G$ be a finite group linearly acting on a polynomial ring $R$ over a field,
such that the characteristic of $R$ does not divide the order of $G$ (\lq\lq non-modular case\rq\rq).
It is well known that the invariant ring $R^G=\{r\in R\co g.r=r\; \forall g\in G\}$ is a finitely
generated sub-algebra of $R$.
In this paper, we provide an algorithm to compute a minimal set of homogeneous
invariant polynomials generating $R^G$. Such generators are also
known as \begriff{fundamental invariants}.

In principal, this can be done as follows: First, one computes primary invariants of $R^G$ and then irreducible secondary invariants. Primary and irreducible secondary invariants together generate
$R^G$ as an algebra, and (potentially after removing some primary invariants) they form an inclusion
minimal generating set~\cite{Kemper}.
N.~Thi\'ery~\cite{Thiery} suggests another algorithm for the computation of a minimal
generating set in the special case of permutation groups, i.e., of groups acting on
$R$ as subgroup of the permutation group of the variables of $R$.
Thi\'ery's algorithm is not based on the computation of primary invariants, but uses
the incremental construction of so-called \emph{SAGBI-Gr\"obner bases}.
His algorithm is implemented in the library \textsc{PerMuVAR} of
\textsc{MuPAD}~\cite{MuPAD}. There is extensive benchmark on \textsc{Magma} and
\textsc{MuPAD}, using the transitive permutation groups on up to nine
variables~\cite{TransPerm}.

Our algorithm comes in one version for permutation groups and one version
for finite matrix groups. We present comparative benchmarks
based on transitive permutation groups on $7$ or $8$ variables.
We implemented our algorithm in a library (i.e., as interpreted code) in
\textsc{Singular}~\cite{Singular}.
In most of the examples, our algorithm is at least $50$ times, often more than
$1000$ times, faster than the
algorithm used by \textsc{Magma}~\cite{magma}.
We also computed minimal generating sets for some transitive permutations groups
on $9$ and $10$ variables. Moreover, we compute minimal generating sets for the natural action of the cyclic groups of order $\le 11$ in characteristic zero and of the cyclic groups
of order $\le 15$ in prime characteristic (but, of course, still in the non-modular case).

We took the key ingredient for our algorithm from a previous paper~\cite{KingSecondary},
where we focused on the computation of secondary invariants of $R^G$. Our algorithm does not involve
solving linear algebra problems that may become rather huge, in contrast to the algorithm exposed
in~\cite{DerksenKemper}. Instead, we use Gr\"obner bases.
Of course, the computation of a Gr\"obner basis can be, in general, a very difficult business.
The main feature of our algorithm is that it involves
at most one computation of a Gr\"obner basis in each degree.
It turns out that this yields a very well-performant algorithm.

Another peculiarity of our algorithm is the fact that it does not rely on \emph{a-priori} bounds
for the maximal degree $\beta(R^G)$ of elements of a minimal generating set of $R^G$. For other
algorithms, like the one presented in~\cite{Thiery}, the performance crucially depends on good
estimates for $\beta(R^G)$. Unfortunately, well known a-priori bounds like Noether's
$\beta(R^G)\le |G|$ are, in general, far from being optimal. In contrast, we rely on more realistic
\emph{a-posteriori} bounds: While incrementally constructing the set of generators, we
obtain informations allowing to estimate $\beta(R^G)$.

We outline our algorithm.
In the case of finite matrix groups, candidates for generators are
found by applying the Reynolds operator to some monomials. In the case of permutation
groups, candidates are found among the \emph{orbit sums}. In increasing degree $d$,
for testing whether a candidate is already contained in the algebra generated by previously found generators,
one computes the normal form with respect to a homogeneous Gr\"obner basis up to degree $d$ of the
ideal spanned by the previously found generators. When starting in a new degree, the Gr\"obner
basis is computed by standard procedures (e.g., Buchberger's algorithm), and when
a new generator of $R^G$ of degree $d$ has been found, one can directly write down a
new Gr\"obner basis up to degree $d$, as we showed in~\cite{KingSecondary}.
Eventually, the ideal spanned by the generators of $R^G$ is $0$--dimensional. Then, $\beta(R^G)$ is
bounded by the highest degree of a monomial not occuring as a leading monomial in
the ideal spanned by the generators. 
Hence, after finishing in that degree, we can stop the quest for more generators.

A modification of our algorithm can be used to compute irreducible secondary invariants. According to our comparative benchmarks, it often performs much better than other known algorithms, including our algorithm presented in~\cite{KingSecondary} and the algorithm recently implemented in \textsc{Magma} V2.13-9 that appears to be not described in a paper yet.

The rest of this paper is organized as follows.
In the next section, we explain our algorithm in more detail. In  Subsection~\ref{sec:benchmarkAlg},
we do some benchmark tests, comparing the implementation of our algorithm in
\textsc{Singular}~\cite{Singular} with the function \texttt{FundamentalInvariants}
of \textsc{Magma}~\cite{magma}. In Subsection~\ref{sec:FurtherCompRes},
we expose some additional examples that seem to be out of reach for other known algorithms.
In the final section, we modify our algorithm in order to compute irreducible secondary invariants, and
do some benchmarks with that algorithm.

\section{The Algorithm}

Let $G$ be a finite group, linearly acting on a polynomial ring $R$ with $n$ variables over some field $K$.
We denote the action of $g\in G$ on $r\in R$ by $g.r\in R$.

Let $R^G=\{r\in R\co g.r=r,\; \forall g\in G\}$ be the invariant ring. Obviously, it is a sub-algebra of $R$,
and we aim at computing a minimal set of generators for $R^G$.
We study here the \begriff{non-modular} case, i.e., the characteristic of $K$ does not divide
the order of $G$. Note that according to~\cite{Kemper}, algorithms for the non-modular case are
useful also in the modular case.

In the non-modular case, we can use the Reynolds operator $\mathrm{Rey}\co R\to R^G$,
that is defined by
\[ \mathrm{Rey}(r) = \frac 1{|G|}\sum_{g\in G} g.r \]
for $r\in R$.
By construction, the restriction of the Reynolds operator to $R^G$ is the identity.
The Reynolds operator does not commute with the ring multiplication. However, it does
commute, if one of the factors is invariant, as in the following lemma. This is, of course, well known. We provide a proof, for completeness.
\begin{lemma}\label{lem:ReynoldsMult}
  Let $p\in R$ and $q\in R^G$. Then, $\mathrm{Rey}(pq) = \mathrm{Rey}(p)q$.
\end{lemma}
\begin{proof}
        For any $g\in G$, we have $g.(pq)=(g.p)(g.q)$. But $q\in R^G$, and thus
        $g.(pq)=(g.p)q$. It follows
        \begin{eqnarray*}
         \mathrm{Rey}(pq) &=& \frac 1{|G|}\sum_{g\in G} g.(pq) \\
         &=& \frac 1{|G|}\sum_{g\in G} (g.p)q = \mathrm{Rey}(p)q
        \end{eqnarray*}
\end{proof}

For any subset $S\subset R$, we denote by $\llangle S\rrangle\subset R$ the sub-algebra generated by $S$,
and by $\langle S\rangle\subset R$ the ideal generated by $S$.
For $d>0$, let $R^G_d$ be the set of homogeneous invariant polynomials of degree $d$.
For an ideal $I\subset R$, let $lm(I)$ be the set of leading monomials occurring in $I$.

For $S\subset R$, let $mon_d(S)\subset R$ be the set of monomials of degree
$d$ that are not contained in $lm(\langle S\rangle)$. This is easy to compute if $S$ is
a homogeneous Gr\"obner basis at least up to degree
$d$.\footnote{The notion of a Gr\"obner basis up to degree $d$ is well known. See, e.g., \cite{KingSecondary} for a definition.}
Let $B_d(S) = \mathrm{Rey}(mon_d(S))$.
By Lemma~3.5.1 and Remark~3.5.3 in~\cite{DerksenKemper}, $B_d(S)$ generates $R^G_d$
as a $K$--vector space.

So, in increasing degree $d$ starting with $d=1$ and $S=\emptyset$,
we may loop through all $b\in B_d(S)$, and add $b$ to the set $S$ of previously
found generators if $b\not\in\langle\langle
S\rangle\rangle$. In that way, one incrementally constructs a generating
set of $R^G$, consisting of homogeneous invariant polynomials. In fact, it
is a \emph{minimal} generating set~\cite{Thiery}.
We can test whether $b \in\langle\langle S\rangle\rangle$ according to the following
lemma. The lemma is well known, but we include a proof for completeness.
\begin{lemma}\label{lem:testAlgebra}
        Let $S\subset R^G$ be a set of homogeneous invariant non-constant polynomials.
        Assume that $R^G_{d'}\subset \langle\langle S\rangle\rangle$ for all $d'<d$, and assume
        that we are in the non-modular case.
        Let $b\in R^G_d$.
        We have $b\in \langle\langle S\rangle\rangle$ if and only if
        $b\in \langle S\rangle$.
\end{lemma}
\begin{proof}
        If $b\in \langle\langle S\rangle\rangle$ then $b\in \langle S\rangle$.
        If $b\in \langle S\rangle$ then we can write $b$ as a finite sum,
        \[ b = \sum_i p_iq_i\]
        with homogeneous polynomials $p_i\in R$ and $q_i\in S$.
        It easily follows from  Lemma~\ref{lem:ReynoldsMult} that
        $b=\mathrm{Rey}(b) = \sum_i \mathrm{Rey}(p_i) q_i$.
        Since the elements of $S$ are non-constant, the $p_i$ are of degree at most $d-1$.
        Hence, $\mathrm{Rey}(p_i)\in R^G_{d'}$ for some $d'<d$. Thus
        $\mathrm{Rey}(p_i)\in \langle \langle S\rangle\rangle$ by hypothesis.
        Therefore, $b\in \langle \langle S\rangle\rangle$.
\end{proof}

As in~\cite{KingSecondary}, we test whether $b\in\langle S\rangle$
by reduction of $b$ with respect to a homogeneous Gr\"obner basis $\mathcal G$
of $\langle S\rangle$ up to degree $d$. After adding $b$ to the set of generators,
we easily obtain a homogeneous Gr\"obner basis up to degree $d$ of
$\langle S\cup\{b\}\rangle$, by the following result from~\cite{KingSecondary}.
Again, we provide its short proof, for completeness.
\begin{theorem}\label{thm:keythm}
  Let $\mathcal G\subset R$ be a homogeneous Gr\"obner basis up to degree $d$ of $\langle \mathcal G\rangle$.
  Let $p\in R$ be a homogeneous polynomial of degree $d$, and $p\not\in \langle \mathcal G\rangle$.
  Then $\mathcal G\cup \{\mathrm{rem}(p;\mathcal G)\}$ is a homogeneous Gr\"obner basis up
  to degree $d$ of $\langle \mathcal G\cup\{p\}\rangle$.
\end{theorem}
\begin{proof}
  Let $r=\mathrm{rem}(p;\mathcal G)$. Since $p\not\in \langle \mathcal G\rangle$ and all polynomials
  are homogeneous, we have $r\not=0$, $\deg(r)=d$, and
  $\langle \mathcal G\cup\{p\}\rangle = \langle \mathcal G\cup
  \{r\}\rangle$.

  By hypothesis, the $S$--polynomials of pairs of elements of $\mathcal G$ are of degree $>d$ or reduce to
  $0$ modulo $\mathcal G$.
  We now consider the $S$--polynomials of $r$ and elements of $\mathcal G$.
  Let $g\in \mathcal G$. 
  By definition of the remainder, we have $lm(r) \not= lm(g)$.
  Therefore the $S$--polynomial of $r$ and $g$ is of degree $>d=\deg(r)$. This implies
  that $\mathcal G\cup \{\mathrm{rem}(p;\mathcal G)\}$ is a homogeneous Gr\"obner basis up
  to degree $d$.
\end{proof}

There is a problem, though. We can incrementally construct a minimal generating
set of $R^G$, in increasing degrees --- but in what degree
shall we stop the construction? By definition, we can stop after having
found the generators in degree $\beta(R^G)$. So, we could adopt a general estimate
for $\beta(R^G)$ like Noether's bound $\beta(R^G)\le |G|$. However, such
general a--priori estimates are very often far from being optimal.

Therefore, we prefer to derive an estimate for $\beta(R^G)$ from the previously
constructed generators. If $S$ is a generating set of $R^G$, then
$\langle S\rangle$ is zero-dimensional, as in the proof of
Proposition~3.3.1 in~\cite{DerksenKemper}.
Hence, there are only finitely many monomials outside $lm(\langle S\rangle)$, of
maximal degree $d_{max}$. Since we can restrict the quest for generators of
$R^G$ of degree $d$ to the Reynolds images of monomials of degree $d$ outside
$lm(\langle S\rangle)$, it follows $\beta(R^G)\le d_{max}$.

Our strategy is to work with a homogeneous Gr\"obner basis of $\langle S\rangle$ that is
subject to a degree restriction, since this is easier to compute than the
entire Gr\"obner basis.
However, for testing whether $\langle S\rangle$ is of dimension $0$, one needs
a Gr\"obner basis of $\langle S\rangle$ \emph{without} degree restriction.
To avoid needless computations, we use the following trick.

By definition, in degree $\beta(R^G)$ we will find a homogeneous
generator of $R^G$, but in degree $\beta(R^G)+1$ we don't.
Hence, only if our incremental construction of $S$ arrives at some degree $d$, such that there
is an element of $S$ in degree $d-2$ but none in degree $d-1$, it makes
sense
to compute a Gr\"obner basis of $\langle S\rangle$ without degree restriction.
If $\dim(\langle S\rangle)=0$, which is tested using the Gr\"obner basis,
then we obtain an estimate for $\beta(R^G)$ that tells us in what degree we can
stop the incremental search.
We thus obtain the following algorithm for the computation of a minimal generating
set of $R^G$, where $G$ is a finite matrix group.

\medskip
\noindent
Algorithm \textsc{Invariant Algebra}
\begin{enumerate}
        \item Construct the Reynolds operator $\mathrm{Rey}\co R\to R^G$.\\
              Let $S=\mathcal G=\emptyset$. Let $d_{max}=0$.
        \item For increasing degree $d$, starting with $d=1$:
        \begin{enumerate}
                \item If $S$ contains elements of degree $d-2$ but no elements of degree $d-1$:
                \begin{enumerate}
                        \item Replace $\mathcal G$ by a (complete) Gr\"obner basis of $\langle S\rangle$.
                        \item If $\dim (\langle S\rangle)=0$ (which is tested using $\mathcal G$),
                                then replace $d_{max}$ by the maximal degree of polynomials
                                outside $lm(\langle S\rangle)$,
                                and if, moreover, $d$ exceeds the new $d_{max}$ then break and return $S$.
                \end{enumerate}
                        If $S$ contains elements of degree $d-1$, replace $\mathcal G$ by
                        a homogeneous Gr\"obner basis $\mathcal G$ of $\langle S\rangle$
                        up to degree $d$.
                \item Compute $B_d(S)$ using $\mathcal G$ and $\mathrm{Rey}$.
                \item For all $b\in B_d(S)$:\\
                 If $\mathrm{rem}(b;\mathcal G)\not=0$
        then replace $S$ by $S\cup \{b\}$ and $\mathcal G$ by
        $\mathcal G\cup \{\mathrm{rem}(b;\mathcal G)\}$.
    \item If $d=d_{max}$ then break and return $S$.
        \end{enumerate}
\end{enumerate}

\noindent
By Theorem~\ref{thm:keythm}, in all steps $\mathcal G$ is a homogeneous Gr\"obner basis of
$\langle S\rangle$ at least up to degree $d$.
Of course, our algorithm has the same basic structure as many other algorithms.
However, our algorithm uses much more elementary methods than the
algorithm described
in~\cite{DerksenKemper} based on linear algebra.
% or in~\cite{Thiery} (based on SAGBI Gr\"obner bases).
No huge systems of linear equations occur, only few explicit Gr\"obner basis computations
are needed (one per degree), and apart from that the most time consuming
operation is the computation of normal forms.
So it is not surprising that usually our implementation of \textsc{Invariant Algebra} in
\textsc{Singular}~\cite{Singular} is much faster than the algorithm from
\cite{DerksenKemper}
implemented in \textsc{Magma}~\cite{magma}.
% or \textsc{MuPAD}~\cite{MuPAD}.

In most of our examples, the computation of homogeneous Gr\"obner bases up to
degree $d$ is not a big deal (there are exceptions, though).
However, for large group orders, the computation of the Reynolds operator exceeds
the ressources. So, the use of the Reynolds operator can be a problem.
In the case of permutation groups, it helps to replace it by
so-called \emph{orbit sums}, which is also used in~\cite{Thiery}. The \emph{orbit}
of a monomial $m\in R$ is $G.m=\{g.m\co g\in G\}$. The \begriff{orbit sum} of $m$
is $m^\circ = \sum_{m'\in G.m} m'$. Of course, $m^\circ\in R^G$.

In contrast to the Reynolds operator, the orbit sums are defined even in the modular case, i.e., if the characteristic of $R$ divides $|G|$.
In the non-modular case, $m^\circ$ is just a scalar multiple of $\mathrm{Rey}(m)$.
In conclusion, if $G$ is a permutation group, we can also define $B_d(S)$ to be the orbit
sums of the monomials in $mon_d(S)$.
Note, however, that even when using orbit sums, the algorithm \textsc{Invariant Algebra}
only works in the non-modular case, since it relies on Lemma~\ref{lem:testAlgebra}.

\section{Computational results}
\label{sec:CompRes}

A classical test bed for the computation of minimal generating sets of invariant
rings of finite groups is provided by \begriff{transitive permuation
groups} \cite{Thiery},~\cite{TransPerm}.
These are groups acting on a polynomial ring $R$ over a field $K$
by permuting variables, such that any two variables are related by the group action.
The \textsc{Magma} function \texttt{TransitveGroups(i)} provides a list of all
classes of transitive permutation groups on $i$ variables.

In our comparative benchmark in Subsection~\ref{sec:benchmarkAlg},
we consider transitive permutation groups on $7$ and $8$
variables in characteristic $0$.
In Subsection~\ref{sec:FurtherCompRes}, we present some more examples of transitive
permutation groups, with up to $11$ variables in characteristic $0$ and up to
$15$ variables in prime characteristic.
Our benchmarks are based on a Linux
x86\underline{\ }64 platform with two AMD Opteron
248 processors (2,2 GHz) and a memory limit of 16~Gb.
%The horizontal
%axis denotes the examples, numbered according to the \textsc{Magma} function
%\texttt{TransitveGroups(i)}. The vertical axis denotes the computation time,
%on a logarithmic scala.

\subsection{Comparative Benchmark based on Transitive Permutation Groups}
\label{sec:benchmarkAlg}

We study here minimal generating sets of invariant rings of transitive permutation groups
on $7$ and $8$ variables, in characteristic $0$.
We compare the following algorithms.
\begin{enumerate}
        \item Our implementation of \textsc{Invariant Algebra} using orbit sums.
          This is part of the \texttt{finvar.lib} library of \textsc{Singular}-3-0-3 (to be
          released soon) and is called
          \texttt{invariant\underline{\ }al\-gebra\underline{\ }perm}.
                We test a $\beta$--version of \textsc{Singular}-3-0-3.
%               The data are given by thick dots in Figure~\ref{fig:benchmark}.
        \item The function \texttt{FundamentalInvariants} of
                \textsc{Magma} V2.13-9 (released January 2007), which, to the best of the author's knowledge, is
                either based on the algorithms described in~\cite{DerksenKemper} or unpublished.
%               The data are given by the crosses in    Figure~\ref{fig:benchmark}.
\end{enumerate}

Note that our implementation in \textsc{Singular} is interpreted code, without any pre-compilation.
As far as known to the author, \texttt{FundamentalIn\-va\-ri\-ants} in
\textsc{Magma} is pre-compiled.

Usually (but not thoroughly) we stopped the computations of an example after two
hours CPU time.
Moreover, we stopped the computation by one algorithm if it took
more than about $1000$ times longer than by the other algorithm. The results are provided
in Table~\ref{tab:MinGenSets7} for the 7 transitive permutation groups on 7 variables,
and in Table~\ref{tab:MinGenSets8} for $45$ transitive permutation groups on 8 variables.
In the first column of the tables, the group is defined by its generators in disjoint cycle presentation.
The rounded CPU times for \textsc{Singular} or \textsc{Magma} in seconds are provided in the next two columns.
The last column of the tables indicates the number of
generators of a minimal generating set of $R^G$, sorted degree-wise.

\begin{table}[hb]
\caption{Transitive permutation groups on 7 variables (characteristic $0$)}
\label{tab:MinGenSets7}
\begin{tabular}{l|r|r|l}
& \textsc{Singular} & \textsc{Magma} & \\
Group & time [s] & time [s] &$\#$ generators (sorted by degree)\\
\hline
$\begin{smallmatrix} (1, 2, 3, 4, 5, 6, 7) \end{smallmatrix}$ & 0.52 & 25.3&{\scriptsize 1,3,8,12,12,6,6}\\
\hline
$\begin{smallmatrix} (1, 2, 3, 4, 5, 6, 7),\\ (1, 6)(2, 5)(3, 4) \end{smallmatrix}$ & 0.67 & 11
&{\scriptsize 1,3,4,6,6,3,3}\\
\hline
$\begin{smallmatrix} (1, 2, 3, 4, 5, 6, 7),\\ (1, 2, 4)(3, 6, 5) \end{smallmatrix}$ & 6.6 & 239
& {\scriptsize 1,1,4,5,8,8,6}\\
\hline
$\begin{smallmatrix} (1, 2, 3, 4, 5, 6, 7),\\ (1, 2)(3, 6) \end{smallmatrix}$ & 16.9 & 107&{\scriptsize 1,1,2,2,2,2,2}\\
\hline
$\begin{smallmatrix} (1, 2, 3, 4, 5, 6, 7),\\ (1, 3, 2, 6, 4, 5) \end{smallmatrix}$ & 81.5 & 600&{\scriptsize 1,1,2,3,4,7,7,5,1}\\
\hline
$\begin{smallmatrix} (3, 4, 5, 6, 7),\\ (1, 2, 3) \end{smallmatrix}$ & 117& 474
&$\begin{smallmatrix} 1,1,1,1,1,1,1,0,0,0, \\ 0,0,0,0,0,0,0,0,0,0,1 \end{smallmatrix}$\\
\hline
$\begin{smallmatrix} (1, 2, 3, 4, 5, 6, 7),\\ (1, 2) \end{smallmatrix}$ & 198 & 0.04
&{\scriptsize 1,1,1,1,1,1,1}\\
\hline

\end{tabular}
\end{table}

\begin{longtable}{l|r|r|l}
        \caption{Transitive permutation groups on 8 variables (characteristic $0$)}
        \label{tab:MinGenSets8}\\
& \textsc{Singular} & \textsc{Magma} & \\
Group & time [s] & time [s] &$\#$ generators (sorted by degree)\\
\hline
\endfirsthead
        Group & \textsc{Singular} & \textsc{Magma} & $\#$ generators (sorted by degree)\\
        \hline
\endhead
        \multicolumn{4}{r}{\emph{Continued on the next page}}
\endfoot

\endlastfoot

$\begin{smallmatrix} (1, 8)(2, 3)(4, 5)(6, 7),\\ (1, 3)(2, 8)(4, 6)(5, 7),\\ (1, 5)(2, 6)(3, 7)(4, 8)\end{smallmatrix}$ & 0.14 & 0.07 & {\scriptsize 1,7,7,7}\\
\hline
$\begin{smallmatrix} (1, 2, 3, 8)(4, 5, 6, 7),\\ (1, 6)(2, 5)(3, 4)(7, 8)\end{smallmatrix}$ & 0.24 & 11.6 & {\scriptsize 1,6,8,12,5}\\
\hline
$\begin{smallmatrix} (1, 2, 3, 8)(4, 5, 6, 7),\\ (1, 5)(2, 6)(3, 7)(4, 8)\end{smallmatrix}$ & 0.35 & 15 & {\scriptsize 1,5,9,16,8}\\
\hline
$\begin{smallmatrix} (1, 8)(2, 3)(4, 5)(6, 7),\\ (1, 3)(2, 8)(4, 6)(5, 7),\\ (1, 5)(2, 6)(3, 7)(4, 8),\\ (4, 5)(6, 7)\end{smallmatrix}$ & 0.35 & 10.8 & {\scriptsize 1,5,5,8,4}\\
\hline
$\begin{smallmatrix} (1, 8)(2, 3)(4, 5)(6, 7),\\ (1, 3)(2, 8)(4, 6)(5, 7),\\ (1, 5)(2, 6)(3, 7)(4, 8),\\ (4, 5)(6, 7),\\ (4, 6)(5, 7)\end{smallmatrix}$ & 0.55 & 34.6 & {\scriptsize 1,4,4,7,3}\\
\hline
$\begin{smallmatrix} (1, 2, 3, 8)(4, 5, 6, 7),\\ (1, 7, 3, 5)(2, 6, 8, 4)\end{smallmatrix}$ & 0.65 & 137 & {\scriptsize 1,4,10,19,15,7}\\
\hline
$\begin{smallmatrix} (1, 8)(2, 3)(4, 5)(6, 7),\\ (1, 3)(2, 8)(4, 6)(5, 7),\\ (1, 5)(2, 6)(3, 7)(4, 8),\\ (2, 3)(4, 5),\\ (2, 3)(6, 7)\end{smallmatrix}$ & 0.65 & 52.2 & {\scriptsize 1,4,4,7,3,1}\\
\hline
$\begin{smallmatrix} (1, 5)(3, 7),\\ (1, 2, 3, 8)(4, 5, 6, 7)\end{smallmatrix}$ & 0.77 & 73.9 & {\scriptsize 1,4,6,11,7,2}\\
\hline
$\begin{smallmatrix} (1, 5)(3, 7),\\ (1, 3, 5, 7)(2, 4, 6, 8),\\ (1, 4, 5, 8)(2, 3, 6, 7)\end{smallmatrix}$ & 0.8 & 167 & {\scriptsize 1,4,6,11,7,3}\\
\hline
$\begin{smallmatrix} (1, 5)(3, 7),\\ (1, 4, 5, 8)(2, 3)(6, 7),\\ (1, 3)(2, 8)(4, 6)(5, 7)\end{smallmatrix}$ & 1.2 & 60.3 & {\scriptsize 1,4,4,6,4,3,2,1}\\
\hline
$\begin{smallmatrix} (4, 8),\\ (1, 8)(2, 3)(4, 5)(6, 7),\\ (1, 3)(2, 8)(4, 6)(5, 7)\end{smallmatrix}$ & 1.4 & 7.38 & {\scriptsize 1,4,4,6,3,1}\\
\hline
$\begin{smallmatrix} (1, 8)(2, 3)(4, 5)(6, 7),\\ (1, 3)(2, 8)(4, 6)(5, 7),\\ (1, 5)(2, 6)(3, 7)(4, 8),\\ (1, 3)(4, 5, 6, 7),\\ (1, 3)(5, 7)\end{smallmatrix}$ & 1.9 & 318 & {\scriptsize 1,3,3,6,3,2,1}\\
\hline
$\begin{smallmatrix} (1, 2, 3, 4, 5, 6, 7, 8)\end{smallmatrix}$ & 2.2 & $>$ 2200 & {\scriptsize 1,4,10,18,16,8,4,4}\\
\hline
$\begin{smallmatrix} (2, 6)(3, 7),\\ (1, 2, 3, 4, 5, 6, 7, 8)\end{smallmatrix}$ & 2.3 & $>$ 2200 & {\scriptsize 1,3,5,8,7,7,4,4}\\
\hline
$\begin{smallmatrix} (2, 6)(3, 7),\\ (1, 2, 3, 8)(4, 5, 6, 7)\end{smallmatrix}$ & 2.3 & 385 & {\scriptsize 1,3,5,9,6,4,2,1}\\
\hline
$\begin{smallmatrix} (1, 8)(2, 3)(4, 5)(6, 7),\\ (1, 3)(2, 8)(4, 6)(5, 7),\\ (1, 5)(2, 6)(3, 7)(4, 8),\\ (1, 3)(4, 5, 6, 7)\end{smallmatrix}$ & 2.4 & 649 & {\scriptsize 1,3,3,7,6,7,5,1}\\
\hline
$\begin{smallmatrix} (1, 2, 3, 4, 5, 6, 7, 8),\\ (1, 5)(3, 7)\end{smallmatrix}$ & 2.8 & $>$ 2800 & {\scriptsize 1,3,7,12,13,9,4,4}\\
\hline
$\begin{smallmatrix} (1, 2, 3, 4, 5, 6, 7, 8),\\ (1, 6)(2, 5)(3, 4)(7, 8)\end{smallmatrix}$ & 3 & 1040 & {\scriptsize 1,4,5,9,8,4,2,2}\\
\hline
$\begin{smallmatrix} (1, 2, 3, 4, 5, 6, 7, 8),\\ (1, 3)(2, 6)(5, 7)\end{smallmatrix}$ & 3.3 & $>$ 3300 & {\scriptsize 1,3,6,11,12,7,2,2}\\
\hline
$\begin{smallmatrix} (4, 8),\\ (1, 2, 3, 8)(4, 5, 6, 7)\end{smallmatrix}$ & 3.7 & 580 & {\scriptsize 1,3,5,8,6,4,2,2}\\
\hline
$\begin{smallmatrix} (2, 6)(3, 7),\\ (1, 3)(5, 7),\\ (1, 2, 3, 4, 5, 6, 7, 8)\end{smallmatrix}$ & 3.7 & $>$ 3600 & {\scriptsize 1,3,3,5,4,4,2,2}\\
\hline
$\begin{smallmatrix} (1, 2, 3, 8),\\ (1, 5)(2, 6)(3, 7)(4, 8)\end{smallmatrix}$ & 4 & $>$ 4000 & {\scriptsize 1,3,4,7,6,4,2,2}\\
\hline
$\begin{smallmatrix} (1, 2, 3, 4, 5, 6, 7, 8),\\ (1, 5)(3, 7),\\ (1, 6)(2, 5)(3, 4)(7, 8)\end{smallmatrix}$ & 4.3 & 5440 & {\scriptsize 1,3,4,7,6,4,2,2}\\
\hline
$\begin{smallmatrix} (1, 8)(2, 3)(4, 5)(6, 7),\\ (1, 3)(2, 8)(4, 6)(5, 7),\\ (1, 5)(2, 6)(3, 7)(4, 8),\\ (1, 2, 3)(4, 6, 5)\end{smallmatrix}$ & 4.9 & 703 & {\scriptsize 1,3,3,7,8,11,7}\\
\hline
$\begin{smallmatrix} (1, 2, 3, 4, 5, 6, 7, 8),\\ (1, 5)(4, 8),\\ (1, 7)(3, 5)(4, 8)\end{smallmatrix}$ & 5 & 4780,6 & {\scriptsize 1,3,3,5,3,3,2,3,1}\\
\hline
$\begin{smallmatrix} (4, 8),\\ (1, 3)(5, 7),\\ (1, 2, 3, 8)(4, 5, 6, 7)\end{smallmatrix}$ & 5.4 & 444 & {\scriptsize 1,3,3,5,3,2,1,1}\\
\hline
$\begin{smallmatrix} (1, 8)(2, 3)(4, 5)(6, 7),\\ (1, 3)(2, 8)(4, 6)(5, 7),\\ (1, 5)(2, 6)(3, 7)(4, 8),\\ (1, 2, 3)(4, 6, 5),\\ (2, 3)(4, 5)\end{smallmatrix}$ & 6.5 & 1995 & {\scriptsize 1,3,3,6,4,3,1}\\
\hline
$\begin{smallmatrix} (2, 6)(3, 7),\\ (1, 3)(4, 8)(5, 7),\\ (1, 2, 3, 8)(4, 5, 6, 7)\end{smallmatrix}$ & 7.5 & $>$ 10800 & {\scriptsize 1,3,3,5,3,2,3,4,3,2,1,1}\\
\hline
$\begin{smallmatrix} (1, 3)(2, 8)(4, 6)(5, 7),\\ (1, 2, 3)(5, 6, 7),\\ (1, 4)(2, 6)(3, 7)(5, 8)\end{smallmatrix}$ & 8.3 & 2410 & {\scriptsize 1,3,3,8,7,9,6,1,1}\\
\hline
$\begin{smallmatrix} (1, 8)(2, 3),\\ (1, 2, 3)(5, 6, 7),\\ (1, 5)(2, 7)(3, 6)(4, 8)\end{smallmatrix}$ & 17.5 & $>$ 7200 & {\scriptsize 1,2,2,5,2,5,4,3,3}\\
\hline
$\begin{smallmatrix} (1, 3)(2, 8),\\ (1, 2, 3),\\ (1, 5)(2, 6)(3, 7)(4, 8)\end{smallmatrix}$ & 31 & $>$ 7200 & {\scriptsize 1,2,2,3,2,3,2,2,1,1}\\
\hline
$\begin{smallmatrix} (1, 8)(2, 3)(4, 5)(6, 7),\\ (1, 3)(2, 8)(4, 6)(5, 7),\\ (1, 5)(2, 6)(3, 7)(4, 8),\\ (1, 2, 3)(4, 6, 5),\\ (2, 5)(3, 4)\end{smallmatrix}$ & 36.5 & $>$ 7200 & {\scriptsize 1,2,2,4,3,5,4,2,2,1,1,1}\\
\hline
$\begin{smallmatrix} (1, 8)(2, 3)(4, 5)(6, 7),\\ (1, 3)(2, 8)(4, 6)(5, 7),\\ (1, 5)(2, 6)(3, 7)(4, 8),\\ (1, 2, 3)(4, 6, 5),\\ (1, 6)(2, 3, 5, 4)\end{smallmatrix}$ & 37 & 3454    & {\scriptsize 1,2,2,4,2,2,1}\\
\hline
$\begin{smallmatrix} (1, 8)(2, 3)(4, 5)(6, 7),\\ (1, 3)(2, 8)(4, 6)(5, 7),\\ (1, 5)(2, 6)(3, 7)(4, 8),\\ (1, 2, 3)(4, 6, 5),\\ (1, 3)(4, 5, 6, 7)\end{smallmatrix}$ & 37 & $>$ 7200 & {\scriptsize 1,2,2,4,2,3,2,2,1}\\
\hline
$\begin{smallmatrix} (1, 3, 5, 7)(2, 4, 6, 8),\\ (1, 3, 8)(4, 5, 7)\end{smallmatrix}$ & 39 & $>$ 7200 & {\scriptsize 1,2,4,8,11,12,7}\\
\hline
$\begin{smallmatrix} (4, 8),\\ (1, 8)(4, 5),\\ (1, 2, 3, 8)(4, 5, 6, 7)\end{smallmatrix}$ & 39 & $>$ 7200 & {\scriptsize 1,2,2,3,2,2,1,1}\\
\hline
$\begin{smallmatrix} (1, 8)(2, 3)(4, 5)(6, 7),\\ (1, 3)(2, 8)(4, 6)(5, 7),\\ (1, 5)(2, 6)(3, 7)(4, 8),\\ (1, 2, 3)(4, 6, 5),\\ (4, 6)(5, 7)\end{smallmatrix}$ & 44 & $>$ 7200 & {\scriptsize 1,2,2,4,3,6,5,5,3}\\
\hline
$\begin{smallmatrix} (1, 3)(2, 8),\\ (1, 2, 3),\\ (1, 8)(4, 5),\\ (1, 5)(2, 6)(3, 7)(4, 8)\end{smallmatrix}$ & 47 & $>$ 7200 & {\scriptsize 1,2,2,3,2,2,1,1,0,0,0,1}\\
\hline
$\begin{smallmatrix} (1, 3)(2, 8),\\ (1, 2, 3),\\ (1, 8)(4, 5),\\ (1, 5)(2, 7, 3, 6)(4, 8)\end{smallmatrix}$ & 50 & $>$ 7200 & {\scriptsize 1,2,2,3,2,2,1,1,0,0,0,0,1,1,1,1}\\
\hline
$\begin{smallmatrix} (4, 8),\\ (1, 8)(2, 3)(4, 5)(6, 7),\\ (1, 2, 3)(5, 6, 7)\end{smallmatrix}$ & 51 & $>$ 7200 & {\scriptsize 1,2,2,3,3,5,4,3,2,1,1,1}\\
\hline
$\begin{smallmatrix} (1, 2, 3, 8),\\ (2, 3),\\ (1, 5)(2, 6)(3, 7)(4, 8)\end{smallmatrix}$ & 51 & 73 & {\scriptsize 1,2,2,3,2,2,1,1}\\
\hline
$\begin{smallmatrix} (1, 5)(4, 8),\\ (1, 8)(2, 3)(4, 5)(6, 7),\\ (1, 2, 3)(5, 6, 7),\\ (2, 3)(4, 8)(6, 7)\end{smallmatrix}$ & 56 & $>$ 7200 & {\scriptsize 1,2,2,3,2,2,1,1,1,3,3,2,2,1,1,1}\\
\hline
$\begin{smallmatrix} (1, 2, 3, 4, 5, 6, 7, 8),\\ (1, 3, 8)(4, 5, 7)\end{smallmatrix}$ & 161.5 & $>$ 7200 & {\scriptsize 1,2,3,5,6,6,5,2}\\
\hline
$\begin{smallmatrix} (1, 2)(3, 4, 5, 6, 7, 8),\\ (1, 2, 3)\end{smallmatrix}$ & 17410 & $>$ 20000 & {$\begin{smallmatrix} 1,1,1,1,1,1,1,1,0,0,0,0,0,0,0,\\
          0,0,0,0,0,0,0,0,0,0,0,0,1\end{smallmatrix}$}\\
\hline
$\begin{smallmatrix} (1, 2, 3, 4, 5, 6, 7, 8),\\ (1, 2) \end{smallmatrix}$ & 24629 & 0.18 & {\scriptsize 1,1,1,1,1,1,1,1}\\
\hline

\end{longtable}

In total, there are $50$ classes of transitive permutation groups on $8$ variables, but
for five of them, neither \textsc{Singular} nor \textsc{Magma} succeeded with the
computation in the realm of our time and memory limits. Note that, according
to~\cite{Thiery}, \textsc{MuPAD} can manage one of these five exceptions with the
library \textsc{PerMuVAR}; with a memory limit of $500$~Mb and a time limit of $2$ days, it can compute $17$ of the $50$ examples.

In the majority of the examples, \textsc{Singular}-3-0-3
is at least $50$ times faster than \textsc{Magma} V2.13-9, in some cases even more
than $1000$ times faster.
There appears to be only one class of exceptions: The symmetric group on $n$ variables
(the last example on Tables~\ref{tab:MinGenSets7} or~\ref{tab:MinGenSets8}, respectively).
This is a special case with a well known theoretical solution.
Since \textsc{Magma} knows that \texttt{TransitiveGroup(7,7)} and
\texttt{TransitiveGroup(8,50)} are symmetric groups, it seems very likely to the
author that \texttt{FundamentalInvariants} simply returns the well known solution in
this case, without computation.
For our algorithm, the invariant ring of the symmetric group is particularly hard, because
the a-posteriori degree bound is not very good. E.g., we find the degree bound $28$ for
the symmetric group on $8$ variables, although a minimal generating set has maximal
degree $8$.

An extensive comparative benchmark of \textsc{MuPAD} and \textsc{Magma}
on transitive permutation groups is provided by~\cite{TransPerm}. There, a different
machine is used, the memory limit is more restrictive (500 Mb), and the time
limit is more generous (2 days).

Note that in the case of small group orders,
it sometimes turned out to be faster to use images of the reynolds operator
(the function \texttt{invariant\underline{\ }algebra\underline{\ }reynolds}
in \textsc{Singular}-3-0-3) rather than orbit sums.
However, for groups of order greater than 1000, 
\textsc{Singular} is hardly able to compute the reynolds operator in reasonable time.
Of course, a pre-compilation would yield a considerable speed-up of our implementation.

\begin{table}[ht]
\caption{Some transitive permutation groups on 9 variables (characteristic $0$)}
\label{tab:MinGenSets9}
\begin{tabular}{l|r|l}
Group & time [s]  & $\#$ generators (sorted by degree)\\
\hline
%Bsp 2:
$\begin{smallmatrix} (1, 2, 9)(3, 4, 5)(6, 7, 8),\\
    (1, 4, 7)(2, 5, 8)(3, 6, 9)
\end{smallmatrix}$ & 6.24 & 1,4,16,24,24 \\
\hline
%Bsp 1:
$\begin{smallmatrix} (1, 2, 3, 4, 5, 6, 7, 8, 9)
\end{smallmatrix}$ & 38.19 & 1,4,14,26,32,18,12,6,6 \\
\hline
%Bsp 5:
$\begin{smallmatrix} (1, 2, 9)(3, 4, 5)(6, 7, 8),\\
    (1, 4, 7)(2, 5, 8)(3, 6, 9),\\
    (1, 2)(3, 6)(4, 8)(5, 7)
\end{smallmatrix}$ & 45.5 & 1,4,8,12,12,10 \\
\hline
%Bsp 4:
$\begin{smallmatrix} (1, 2, 9)(3, 4, 5)(6, 7, 8),\\
    (1, 2)(4, 5)(7, 8),\\
    (1, 4, 7)(2, 5, 8)(3, 6, 9)
\end{smallmatrix}$ & 55.3 & 1,3,10,14,19,9,2 \\
\hline
%Bsp 7:
$\begin{smallmatrix} (1, 2, 9)(3, 4, 5)(6, 7, 8),\\
    (1, 4, 7)(2, 5, 8)(3, 6, 9),\\
    (3, 4, 5)(6, 8, 7)
\end{smallmatrix}$ & 84.3 & 1,2,8,9,16,18,14,4,2 \\
\hline
%Bsp 3:
$\begin{smallmatrix} (1, 2, 3, 4, 5, 6, 7, 8, 9),\\
    (1, 8)(2, 7)(3, 6)(4, 5)
\end{smallmatrix}$ & 141.6 & 1,4,7,13,16,12,6,3,3 \\
\hline
%Bsp 8:
$\begin{smallmatrix} (1, 2, 9)(3, 4, 5)(6, 7, 8),\\
    (1, 2)(4, 5)(7, 8),\\
    (1, 4, 7)(2, 5, 8)(3, 6, 9),\\
    (3, 6)(4, 7)(5, 8)
\end{smallmatrix}$ & 280.7 & 1,3,6,8,9,8,2 \\
\hline
%Bsp 12:
$\begin{smallmatrix} (3, 4, 5)(6, 8, 7),\\
    (1, 4, 7)(2, 5, 8)(3, 6, 9),\\
    (3, 6)(4, 7)(5, 8)
\end{smallmatrix}$ & 290.5 & 1,2,6,6,9,8,4 \\
\hline
%Bsp 6:
$\begin{smallmatrix} (1, 4, 7)(2, 8, 5),\\
    (1, 2, 3, 4, 5, 6, 7, 8, 9)
\end{smallmatrix}$ & 455.1 & 1,2,6,11,20,25,26,10,8 \\
\hline
\end{tabular}
\end{table}

\begin{table}[ht]
\caption{Some transitive permutation groups on 10 variables (characteristic $0$)}
\label{tab:MinGenSets10}
\begin{tabular}{l|r|l}
Group &  time [s] & $\#$ generators (sorted by degree)\\
\hline
%Bsp 2:
$\begin{smallmatrix} (1, 3, 5, 7, 9)(2, 4, 6, 8, 10),\\
    (1, 4)(2, 3)(5, 10)(6, 9)(7, 8)
\end{smallmatrix}$ & 12.3 &  1,7,14,29,28,12\\
\hline
%Bsp 1:
$\begin{smallmatrix} (1, 2, 3, 4, 5, 6, 7, 8, 9, 10)
\end{smallmatrix}$ & 306 &  1,5,16,36,48,32,12,8,4,4\\
\hline
%Bsp 8:
$\begin{smallmatrix} (2, 7)(5, 10),\\
    (1, 3, 5, 7, 9)(2, 4, 6, 8, 10)
\end{smallmatrix}$ & 478 &  1,3,8,14,21,16,12,8,4,3\\
\hline
%Bsp 4:
$\begin{smallmatrix} (1, 3, 5, 7, 9)(2, 4, 6, 8, 10),\\
    (1, 2, 9, 8)(3, 6, 7, 4)(5, 10)
\end{smallmatrix}$ & 1294 &  1,4,9,20,31,23,8\\
\hline
%Bsp 3:
$\begin{smallmatrix} (1, 2, 3, 4, 5, 6, 7, 8, 9, 10),\\
    (1, 8)(2, 7)(3, 6)(4, 5)(9, 10)
\end{smallmatrix}$ & 1425 &  1,5,8,18,24,17,6,4,2,2\\
\hline
%Bsp xx:
%$\begin{smallmatrix}
%\end{smallmatrix}$ &  &  \\
%\hline
\end{tabular}
\end{table}

\subsection{Further computational results}
\label{sec:FurtherCompRes}

In this subsection, we consider some more examples of transitive
permutation groups, acting
on up to 15 variables.
Given the results exposed in the preceding subsection, it seems very unlikely to us
that \textsc{Magma} V2.13-9 is able to compute these examples in reasonable time. Hence, we only tried with \textsc{Singular}-3-0-3 (Beta version).
Table~\ref{tab:MinGenSets9} and Table~\ref{tab:MinGenSets10} provide the results for some transitive permutation groups on 9 and 10 variables, in characteristic $0$; here, we used orbit sums. According to~\cite{Thiery}, \textsc{MuPAD} can handle $5$ of the
transitive permutation groups on $9$ variables (in total, there are $34$ of them) using the library \textsc{PerMuVAR}, with a memory limit of $500$~Mb and a time limit of $2$ days.

\begin{table}[ht]
\caption{Natural action of $C_n$ on $n$ variables (characteristic $0$)}
\label{tab:cyclic}
\begin{tabular}{l|r|r|l}
$n$ &  time [s] & mem. [Mb] & $\#$ generators (sorted by degree)\\
\hline
6 & 0.05 & 0.746& 1,3,6,6,2,2\\
7 & 0.17 & 1.25 & 1,3,8,12,12,6,6\\
8 & 1.54 & 2.25 & 1,4,10,18,16,8,4,4\\
9 & 35.6 & 11.92& 1,4,14,26,32,18,12,6,6\\
10& 298.3& 54.16& 1,5,16,36,48,32,12,8,4,4\\
11& 1187 & 116  & 1,5,20,50,82,70,50,30,20,10,10\\
\end{tabular}
\end{table}

\begin{table}[ht]
\caption{Natural action of $C_n$ on $n$ variables (characteristic $p>0$)}
\label{tab:cyclicmodp}
\begin{tabular}{l|l|r|r|l}
$n$ & $p$ & time [s] & mem. [Mb] & $\#$ generators (sorted by degree)\\
\hline
6 & 5 & 0.03 & 0.746& 1,3,6,6,2,2\\
7 & 2 & 0.09 & 0.746& 1,3,8,12,12,6,6\\
8 & 3 & 0.34 & 1.25 & 1,4,10,18,16,8,4,4\\
9 & 2 & 1.65 & 1.86 & 1,4,14,26,32,18,12,6,6\\
10& 3 & 12.7 & 4.48 & 1,5,16,36,48,32,12,8,4,4\\
11& 2 & 73.5 & 9.33 & 1,5,20,50,82,70,50,30,20,10,10\\
12& 5 & 693  & 33.2 & 1,6,24,64,104,84,36,20,12,8,4,4\\
13& 2 & 4079 & 81.1 & 1,6,28,84,168,180,132,84,60,36,24,12,12\\
14 & 3& 25280& 304.3& 1,7,32,104,216,242,162,96,42,30,18,12,6,6\\
15 & 2& 99873& 780.4& 1,7,38,130,306,388,264,120,88,56,40,24,16,8,8
\end{tabular}
\end{table}

A rather harmlessly looking class of transitive permutation groups is the natural action
of the cyclic group $C_n$ of order $n$ on $n$ variables. The maximal degree occuring
in a minimal generating set is, by Noether's bound, of course at most $|C_n|=n$, hence, quite small.
However, the minimal number of generators of $R^{C_n}$ is surprisingly large.
Since here the group orders are very small,
we use the Reynolds operator rather than orbit sums for the generation of invariants.
For $n\le 5$ the computation is finished in almost no time, so we omit them in our tables. Table~\ref{tab:cyclic} provides the result for $n=6,...,11$ in characteristic $0$.
Recall that for the timings in Tables~\ref{tab:MinGenSets7}--\ref{tab:MinGenSets10} we
used orbit sums and not the Reynolds operator --- this explains the different computation
times in the case of cyclic groups.

Table~\ref{tab:cyclicmodp} provides the results for $n=6,...,15$ in small prime
characteristic $p>0$, of course such that $p$ does not divide $n$ (non-modular case).
Apparently this is much easier than characteristic $0$. The reason is that in characteristic
$0$ the coefficients occuring in the Gr\"obner bases become very huge. By consequence, it takes too long to compute normal forms.

Note that the in all examples, the number of generators in each degree is the same in
characteristic 0 and in non-modular prime characteristic. It is in fact conjectured
that this is always the case~\cite{ThieryPersonal}.

To work in prime characteristic is not the only way to simplify the computations.
As a last example, we study here the action of $S_5$ on pairs, which yields a $10$--dimensional
representation of $S_5$. One can decompose it into a $1$-, a $4$- and a $5$--dimensional
irreducible representation, and in this form, the representation is given by
the matrices \label{expl:S5D10}
   \begin{eqnarray*}
    M_1 &=& \left(\begin{smallmatrix}
1 &  0  & 0  & 0  & 0  & 0  & 0 &  0 &  0 &  0\\
0 &  1 &\frac 13 &\frac 13 &\frac 13&   0 &  0&   0&   0&   0\\
0 &  0 &\frac 13&-\frac 23&-\frac 23 &  0&   0 &  0&   0&   0\\
0 &  0&-\frac 23& \frac 13&-\frac 23&   0&   0&   0&   0&   0\\
0 &  0&-\frac 23&-\frac 23& \frac 13&   0&   0&   0&   0&   0\\
0 &  0&   0&   0&   0&   1&   0&   0&   0&   0\\
0 &  0&   0&   0 &  0 &  0 &  0 &  0&   1 &  0\\
0 &  0&   0 &  0 &  0 &  0 &  0 &  0 &  0 &  1\\
0 &  0 &  0 &  0 &  0 &  0 &  1 &  0 &  0 &  0\\
0 &  0 &  0 &  0 &  0 &  0  & 0 &  1 &  0 &  0
              \end{smallmatrix}\right) \\
    M_2 &=& \left(\begin{smallmatrix}
1 &  0 &  0 &  0 &  0 &  0 &  0&   0 &  0  & 0\\
0 &  0 &\frac 13&-\frac 23&-\frac 23&   0&   0&  0 &  0 &  0\\
0 &  0&-\frac 23& \frac 13&-\frac 23&   0 &  0&  0 &  0 &  0\\
0 &  0&-\frac 23&-\frac 23& \frac 13&   0 &  0&   0 &  0 &  0\\
0 &  1& \frac 13& \frac 13& \frac 13&   0 &  0&   0 &  0 &  0\\
0 &  0 &  0 &  0 &  0&  -1&  -1&   1&   1 &  0\\
0 &  0 &  0 &  0 &  0 & -1 &  0&   0&   0 &  1\\
0 &  0 &  0 &  0&   0 & -1&   0 &  1 &  0  & 0\\
0 &  0 &  0 &  0 &  0 &  0&  -1&  0  & 0 &  0\\
0  & 0 &  0 &  0 &  0 &  0&  -1&  1 &  0 &  0
              \end{smallmatrix}\right)
   \end{eqnarray*}
According to an advice of G. Kemper~\cite{KemperPersonal}, we used this as an example for the computation of irreducible secondary invariants (see~\cite{KingSecondary} or
the benchmark in the next section). But of course it is also a nice example for the
computation of a minimal generating set.

We could describe that representation of $S_5$ by a transitive permutation group
on $10$ variables. However, in that formulation of the problem, our algorithm would take a very long time to find a minimal generating set.
But after the decomposition, our algorithm \textsc{Invariant Algebra}
executed in \textsc{Singular} 3-0-2 finds a minimal generating set after 
$47.8$ minutes using $4.4$~Gb in characteristic $0$ respectively after 
only $84.2$ seconds using $81.7$~Mb
in characteristic $7$. In both cases, there is a minimal number of
$1,2,4,7,10,13,13,4,2$ generators sorted by degree.

Even using the decomposition, the \textsc{Magma} V2.13-9 function 
\texttt{FundamentalInvariants} is unable to find a minimal generating set 
in less than $4$ hours, both in characteristic $0$ and in characteristic 
$7$.

\section{Application to irreducible secondary invariants}

In~\cite{KingSecondary}, we presented an algorithm
for the computation of secondary invariants and a specialised version for the computation
of \emph{irreducible} secondary invariants. Shortly after the first version of~\cite{KingSecondary} was posted, there was a new release of \textsc{Magma} containing
a new algorithm of G. Kemper for the computation of secondary invariants.
Unfortunately, to the best of the author's knowledge, Kemper did not describe his
new algorithm in a manuscript, yet. So it is not clear how it differs from the
algorithm described in~\cite{Kemper}, \cite{KemperSteel} and~\cite{DerksenKemper}
or from the algorithm described in~\cite{KingSecondary}.

By a slight modification, our algorithm can be used to compute irreducible secondary invariants.
For this, let $P$ be a system of primary invariants. In Step~(1) of algorithm
\textsc{Invariant Algebra}, let $S=P$ and let $G$ be a Gr\"obner basis of $P$.
The rest of the algorithm remains unchanged. In the end, it returns the
union of $P$
with a system of irreducible secondary invariants.
Note that this algorithm does not involve an application of Molien's Theorem. So, it applies
also to cases when the Molien series is difficult to compute.

In the new version
of \texttt{irred\underline{\ }secondary\underline{\ }char0} in \textsc{Singular}-3-0-3,
we combine both algorithms, i.e., we use the Molien series and power products as described in~\cite{KingSecondary} in low degrees, and the algorithm \textsc{Invariant Algebra} in higher degrees.

For our benchmark, we use Expl.~(4)--(9) from~\cite{KingSecondary},
and one additional example, that appeared
in our study of ideal Turaev--Viro invariants (see~\cite{KingIdeal2} or~\cite{KingIdeal1}
for background material).
Expl.~(9) is the $10$--dimensional representation of $S_5$ discussed above on
Page~\pageref{expl:S5D10}; primary invariants can be easily found by considering
the direct summands of the decomposition separately.
For the sake of brevity, we do not re-define the other examples
from~\cite{KingSecondary}, but just provide the new example.
The ring variables are called $x_1,x_2,...$.
Let $e_i$ be the column vector with $1$ in position $i$ and $0$ otherwise.
In all examples of this section, we work in characteristic $0$.

\begin{itemize}
        \item[(10)] A $20$--dimensional representation of $S_3$ is given by the
        matrices \nopagebreak
   \begin{eqnarray*}
    M_1 &=& \left(
    e_{2}e_{1}e_{3}e_{19}e_{9}e_{13}e_{17}e_{11}e_{5}e_{15}
    e_{8}e_{16}e_{6}e_{14}e_{10}e_{12}e_{7}e_{20}e_{4}e_{18}
\right)\\
% 0,1,0,0,0,0,0,0,0,0,0,0,0,0,0,0,0,0,0,0,  1
% 1,0,0,0,0,0,0,0,0,0,0,0,0,0,0,0,0,0,0,0,
% 0,0,1,0,0,0,0,0,0,0,0,0,0,0,0,0,0,0,0,0,
% 0,0,0,0,0,0,0,0,0,0,0,0,0,0,0,0,0,0,1,0,
% 0,0,0,0,0,0,0,0,1,0,0,0,0,0,0,0,0,0,0,0,  5
% 0,0,0,0,0,0,0,0,0,0,0,0,1,0,0,0,0,0,0,0,
% 0,0,0,0,0,0,0,0,0,0,0,0,0,0,0,0,1,0,0,0,
% 0,0,0,0,0,0,0,0,0,0,1,0,0,0,0,0,0,0,0,0,
% 0,0,0,0,1,0,0,0,0,0,0,0,0,0,0,0,0,0,0,0,
% 0,0,0,0,0,0,0,0,0,0,0,0,0,0,1,0,0,0,0,0, 10
% 0,0,0,0,0,0,0,1,0,0,0,0,0,0,0,0,0,0,0,0,
% 0,0,0,0,0,0,0,0,0,0,0,0,0,0,0,1,0,0,0,0,
% 0,0,0,0,0,1,0,0,0,0,0,0,0,0,0,0,0,0,0,0,
% 0,0,0,0,0,0,0,0,0,0,0,0,0,1,0,0,0,0,0,0,
% 0,0,0,0,0,0,0,0,0,1,0,0,0,0,0,0,0,0,0,0, 15
% 0,0,0,0,0,0,0,0,0,0,0,1,0,0,0,0,0,0,0,0,
% 0,0,0,0,0,0,1,0,0,0,0,0,0,0,0,0,0,0,0,0,
% 0,0,0,0,0,0,0,0,0,0,0,0,0,0,0,0,0,0,0,1,
% 0,0,0,1,0,0,0,0,0,0,0,0,0,0,0,0,0,0,0,0,
% 0,0,0,0,0,0,0,0,0,0,0,0,0,0,0,0,0,1,0,0; 20
                M_2&=&\left(
                e_{1}e_{3}e_{2}e_{4}e_{6}e_{5}e_{10}e_{9}e_{8}e_{7}
                e_{13}e_{16}e_{11}e_{19}e_{20}e_{12}e_{18}e_{17}e_{14}e_{15}
                \right)
% 1,0,0,0,0,0,0,0,0,0,0,0,0,0,0,0,0,0,0,0,  1
% 0,0,1,0,0,0,0,0,0,0,0,0,0,0,0,0,0,0,0,0,
% 0,1,0,0,0,0,0,0,0,0,0,0,0,0,0,0,0,0,0,0,
% 0,0,0,1,0,0,0,0,0,0,0,0,0,0,0,0,0,0,0,0,
% 0,0,0,0,0,1,0,0,0,0,0,0,0,0,0,0,0,0,0,0,  5
% 0,0,0,0,1,0,0,0,0,0,0,0,0,0,0,0,0,0,0,0,
% 0,0,0,0,0,0,0,0,0,1,0,0,0,0,0,0,0,0,0,0,
% 0,0,0,0,0,0,0,0,1,0,0,0,0,0,0,0,0,0,0,0,
% 0,0,0,0,0,0,0,1,0,0,0,0,0,0,0,0,0,0,0,0,
% 0,0,0,0,0,0,1,0,0,0,0,0,0,0,0,0,0,0,0,0, 10
% 0,0,0,0,0,0,0,0,0,0,0,0,1,0,0,0,0,0,0,0,
% 0,0,0,0,0,0,0,0,0,0,0,0,0,0,0,1,0,0,0,0,
% 0,0,0,0,0,0,0,0,0,0,1,0,0,0,0,0,0,0,0,0,
% 0,0,0,0,0,0,0,0,0,0,0,0,0,0,0,0,0,0,1,0,
% 0,0,0,0,0,0,0,0,0,0,0,0,0,0,0,0,0,0,0,1, 15
% 0,0,0,0,0,0,0,0,0,0,0,1,0,0,0,0,0,0,0,0,
% 0,0,0,0,0,0,0,0,0,0,0,0,0,0,0,0,0,1,0,0,
% 0,0,0,0,0,0,0,0,0,0,0,0,0,0,0,0,1,0,0,0,
% 0,0,0,0,0,0,0,0,0,0,0,0,0,1,0,0,0,0,0,0,
% 0,0,0,0,0,0,0,0,0,0,0,0,0,0,1,0,0,0,0,0; 20
\end{eqnarray*}
\nopagebreak
        We use the following sub-optimal primary invariants:
        \allowdisplaybreaks
        \begin{eqnarray*}\allowdisplaybreaks
&&x_{1}+x_{2}+x_{3},\;\;x_{1}x_{2}+x_{1}x_{3}+x_{2}x_{3},\;\;x_{1}x_{2}x_{3},\;\;x_{4}+x_{14}+x_{19},\\
&&x_{4}x_{14}+x_{4}x_{19}+x_{14}x_{19},\;\;x_{4}x_{14}x_{19},\;\;x_{5}+x_{6}+x_{8}+x_{9}+x_{11}+x_{13},\\
&&x_{8}x_{9}+x_{5}x_{11}+x_{6}x_{13},\;\;x_{6}x_{8}+x_{5}x_{9}+x_{11}x_{13},\\
&&x_{5}x_{8}+x_{6}x_{9}+x_{6}x_{11}+x_{9}x_{11}+x_{5}x_{13}+x_{8}x_{13},\\
&&x_{5}x_{6}x_{11}+x_{5}x_{8}x_{11}+x_{8}x_{9}x_{11}+x_{5}x_{6}x_{13}+x_{6}x_{9}x_{13}+x_{8}x_{9}x_{13},\\
&&x_{5}^6+x_{6}^6+x_{8}^6+x_{9}^6+x_{11}^6+x_{13}^6,\;\;x_{12}+x_{16},\;\;x_{12}x_{16},\\
&&x_{7}+x_{10}+x_{15}+x_{17}+x_{18}+x_{20},\;\;x_{7}x_{17}+x_{10}x_{18}+x_{15}x_{20},\\
&&x_{10}x_{15}+x_{17}x_{18}+x_{7}x_{20},\;\;
x_{7}x_{15}+x_{10}x_{17}+x_{7}x_{18}+x_{15}x_{18}+x_{10}x_{20}+x_{17}x_{20},\\
&&x_{7}x_{10}x_{17}+x_{7}x_{15}x_{17}+x_{7}x_{10}x_{18}+x_{15}x_{17}x_{20}+x_{10}x_{18}x_{20}+x_{15}x_{18}x_{20},\\
&&x_{7}^6+x_{10}^6+x_{15}^6+x_{17}^6+x_{18}^6+x_{20}^6
        \end{eqnarray*}
\end{itemize}

In this example, there are $248832$ secondary invariants of maximal degree $26$,
among wich are $283$ irreducible secondary invariants of maximal degree $4$.
The sheer number of secondary invariants (which can be computed by Molien's Theorem) makes the computations hardly manageable for any algorithm that is based on the generation of power products, as the one described in~\cite{Kemper},
\cite{KemperSteel} and~\cite{DerksenKemper}, or the one described in~\cite{KingSecondary}.
It is in fact too much for \textsc{Magma}~V2.13-9 and for \textsc{Singular}-3-0-2.
However, our new algorithm implemented in \textsc{Singular}-3-0-3 just needs few seconds
to compute all irreducible secondary invariants.

In Table~\ref{tab:NewBench}, we compare a $\beta$--version of \textsc{Singular}-3-0-3
(function \texttt{irred\underline{\ }se\-con\-da\-ry\underline{\ }char0}) with \textsc{Magma}~V2.13-9 (function \texttt{Irredu\-cible\-Secon\-dary\-Inva\-riants},
released in January, 2007).
The only exception is Example~(9), that we compute with our new algorithm, but 
based on \textsc{Singular}-3-0-2.
For convenience, we repeat in Table~\ref{tab:NewBench} the timings for \textsc{Singular}-3-0-2 and \textsc{Magma}~V2.13-8 from~\cite{KingSecondary}.

\begin{table}[h]
\caption{Comparative benchmark for the computation of irreducible secondary invariants}
\label{tab:NewBench}
\begin{tabular}{l|c|c|c|c|}
& \textsc{Singular}  &\textsc{Magma}& \textsc{Magma}&\textsc{Singular}  \\
&3-0-3$\beta$ &V2.13-9 & V2.13-8  & 3-0-2   \\

\hline  
Expl. (4) & 0.07 s    & 0.09        &0.48 s     &0.32 s   \\
          & 0.91 Mb   & 7.35 Mb     &9.09 Mb    &2.97 Mb  \\
\hline
Expl. (5) & 7.75 s    & 0.49 s      &6.66 s     &9.69 s   \\
          & 10.9 Mb   & 9.06 Mb     &31.82 Mb   &17.0 Mb  \\
\hline
Expl. (6) & 1.63 s    & 2.49 s      &118.51 s   &16.55 s  \\
          & 6.9 Mb    & 19.8 Mb     &54.0 Mb    &39.0 Mb  \\
\hline
Expl. (7) & 0.34 s    & 36.57 s     &$> 7$ h    &20.94 s  \\
          & 2.52 Mb   & 30.1 Mb     &$> 15$ Gb  &35.1 Mb  \\
\hline
Expl. (8) & 1.05 s    & $>$ 72 min  & ---       & 50.7 min\\
          & 7.08 Mb   & $>$ 2.5 Gb  &(259.5 Gb) & 3.36 Gb \\
\hline
Expl. (9)* & 17.2 min  & 29.9 min    & ---       & 99.2 min\\
               &  4.67 Gb  & 399.5 Mb    & (55.62 Gb)& 7.35 Gb \\
\hline
Expl. (10)& 6.83 s    & $>$ 280 min &           & ---     \\
          & 29.4 Mb   & $>$ 9.9 Gb  &           & ---     \\
\hline
\end{tabular}
\end{table}

The outcome of these benchmarks is less clear than of our benchmarks on
minimal generating sets.
In $3$ of the 7 examples, our algorithm and the one used in \textsc{Magma} V2.13-9
show more or less the same performance (by factors less than $2$),
in one example \textsc{Magma} is faster by a factor of about $15$,
whereas in $3$ examples our algorithm is faster by factors
between 100 and at least 4000.

Note that  in Expl.~(9), the critical part is the computation of a Gr\"obner basis
of primary and irreducible secondary invariants. The rest of the computations just
takes about 5 minutes. 
The beta version of \textsc{Singular}-3-0-3 spends much more than 
$30$ minutes with the computation of a Gr\"obner basis. Here, the old version 
\textsc{Singular}-3-0-2 happens to be quicker.

\end{document}